\documentclass[12pt]{article}
\usepackage{bm,subeqnarray,eqsection,indent,cite}
\usepackage{latexsym}
\usepackage{amsmath}           
\usepackage{amsfonts}          
\usepackage{amssymb}           

\begin{document}

\title{When is a Riesz distribution a complex measure?}

\author{
  {\small Alan D.~Sokal\thanks{Also at Department of Mathematics,
           University College London, London WC1E 6BT, England.}}  \\[-1.7mm]
  {\small\it Department of Physics}       \\[-1.7mm]
  {\small\it New York University}         \\[-1.7mm]
  {\small\it 4 Washington Place}          \\[-1.7mm]
  {\small\it New York, NY 10003 USA}      \\[-1.7mm]
  {\small\tt sokal@nyu.edu}           \\[-1.7mm]
  {\protect\makebox[5in]{\quad}}  
}

\date{June 12, 2009 \\[1mm] revised May 20, 2010}
\maketitle
\thispagestyle{empty}   

\vspace{-5mm}

\begin{abstract}
Let $\mathcal{R}_\alpha$ be the Riesz distribution on a
simple Euclidean Jordan algebra, parametrized by $\alpha \in \mathbb C$.
I give an elementary proof of the necessary and sufficient condition
for $\mathcal{R}_\alpha$
to be a locally finite complex measure (= complex Radon measure).

\medskip

Soit $\mathcal{R}_\alpha$ la distribution de Riesz sur une
alg\`ebre de Jordan euclidienne simple,
param\'etris\'ee par $\alpha \in \mathbb C$.
Je donne une d\'emonstration \'el\'ementaire
de la condition n\'ecessaire et suffisante
pour que $\mathcal{R}_\alpha$ soit une mesure complexe localement finie
(= mesure de Radon complexe).
\end{abstract}


\bigskip
\noindent
{\bf Mathematics Subject Classification (MSC 2000) codes:}
43A85 (Primary);
17A15, 17C99, 28C10, 44A10, 46F10, 47G10, 60E05, 62H05 (Secondary).

\clearpage

\newtheorem{defin}{Definition}[section]
\newtheorem{definition}[defin]{Definition}
\newtheorem{prop}[defin]{Proposition}
\newtheorem{proposition}[defin]{Proposition}
\newtheorem{lem}[defin]{Lemma}
\newtheorem{lemma}[defin]{Lemma}
\newtheorem{guess}[defin]{Conjecture}
\newtheorem{ques}[defin]{Question}
\newtheorem{question}[defin]{Question}
\newtheorem{prob}[defin]{Problem}
\newtheorem{problem}[defin]{Problem}
\newtheorem{thm}[defin]{Theorem}
\newtheorem{theorem}[defin]{Theorem}
\newtheorem{cor}[defin]{Corollary}
\newtheorem{corollary}[defin]{Corollary}
\newtheorem{conj}[defin]{Conjecture}
\newtheorem{conjecture}[defin]{Conjecture}
\newtheorem{examp}[defin]{Example}
\newtheorem{example}[defin]{Example}
\newtheorem{claim}[defin]{Claim}

\renewcommand{\theenumi}{\alph{enumi}}
\renewcommand{\labelenumi}{(\theenumi)}
\def\prf{\par\noindent{\bf Proof.\enspace}\rm}
\def\rmk{\par\medskip\noindent{\bf Remark.\enspace}\rm}

\newcommand{\be}{\begin{equation}}
\newcommand{\ee}{\end{equation}}
\newcommand{\<}{\langle}
\renewcommand{\>}{\rangle}
\newcommand{\widebar}{\overline}
\def\reff#1{(\protect\ref{#1})}
\def\spose#1{\hbox to 0pt{#1\hss}}
\def\ltapprox{\mathrel{\spose{\lower 3pt\hbox{$\mathchar"218$}}
 \raise 2.0pt\hbox{$\mathchar"13C$}}}
\def\gtapprox{\mathrel{\spose{\lower 3pt\hbox{$\mathchar"218$}}
 \raise 2.0pt\hbox{$\mathchar"13E$}}}
\def\textprime{${}^\prime$}
\def\proof{\par\medskip\noindent{\sc Proof.\ }}
\newcommand{\qed}{\quad $\Box$ \medskip \medskip}
\def\firstproof{\par\medskip\noindent{\sc First proof.\ }}
\def\secondproof{\par\medskip\noindent{\sc Second proof.\ }}
\def\proofof#1{\bigskip\noindent{\sc Proof of #1.\ }}
\def\secondproofof#1{\bigskip\noindent{\sc Second proof of #1.\ }}
\def\alternateproofof#1{\bigskip\noindent{\sc Alternate proof of #1.\ }}
\def\half{ {1 \over 2} }
\def\third{ {1 \over 3} }
\def\twothird{ {2 \over 3} }
\def\smfrac#1#2{{\textstyle{#1\over #2}}}
\def\smsmfrac#1#2{{\scriptstyle{#1\over #2}}}
\def\smhalf{ {\smfrac{1}{2}} }
\def\smsmhalf{ {\smsmfrac{1}{2}} }
\newcommand{\real}{\mathop{\rm Re}\nolimits}
\renewcommand{\Re}{\mathop{\rm Re}\nolimits}
\newcommand{\imag}{\mathop{\rm Im}\nolimits}
\renewcommand{\Im}{\mathop{\rm Im}\nolimits}
\newcommand{\sgn}{\mathop{\rm sgn}\nolimits}
\newcommand\supp{\mathop{\rm supp}\nolimits}
\newcommand{\diag}{\mathop{\rm diag}\nolimits}
\newcommand{\pf}{\mathop{\rm pf}\nolimits}
\newcommand{\hf}{\mathop{\rm hf}\nolimits}
\newcommand{\tr}{\mathop{\rm tr}\nolimits}
\newcommand{\Det}{\mathop{\rm Det}\nolimits}
\newcommand{\per}{\mathop{\rm per}\nolimits}
\newcommand{\adj}{\mathop{\rm adj}\nolimits}
\newcommand{\Res}{\mathop{\rm Res}\nolimits}
\newcommand{\Jtr}{\mathop{\rm Jtr}\nolimits}
\newcommand{\Jdet}{\mathop{\rm Jdet}\nolimits}
\newcommand{\Mdet}{\mathop{\rm Mdet}\nolimits}
\def\hboxscript#1{ {\hbox{\scriptsize\em #1}} }
\def\hboxrm#1{ {\hbox{\scriptsize\rm #1}} }

\newcommand{\restrict}{\upharpoonright}
\renewcommand{\emptyset}{\varnothing}

\def\Z{{\mathbb Z}}
\def\ZZ{{\mathbb Z}}
\def\R{{\mathbb R}}
\def\C{{\mathbb C}}
\def\CC{{\mathbb C}}
\def\HH{{\mathbb H}}
\def\OO{{\mathbb O}}
\def\N{{\mathbb N}}
\def\NN{{\mathbb N}}
\def\Q{{\mathbb Q}}

\newcommand{\scra}{{\mathcal{A}}}
\newcommand{\scrb}{{\mathcal{B}}}
\newcommand{\scrc}{{\mathcal{C}}}
\newcommand{\scrd}{{\mathcal{D}}}
\newcommand{\scre}{{\mathcal{E}}}
\newcommand{\scrf}{{\mathcal{F}}}
\newcommand{\scrg}{{\mathcal{G}}}
\newcommand{\scrh}{{\mathcal{H}}}
\newcommand{\scri}{{\mathcal{I}}}
\newcommand{\scrj}{{\mathcal{J}}}
\newcommand{\scrk}{{\mathcal{K}}}
\newcommand{\scrl}{{\mathcal{L}}}
\newcommand{\scrm}{{\mathcal{M}}}
\newcommand{\scrn}{{\mathcal{N}}}
\newcommand{\scro}{{\mathcal{O}}}
\newcommand{\scrp}{{\mathcal{P}}}
\newcommand{\scrq}{{\mathcal{Q}}}
\newcommand{\scrr}{{\mathcal{R}}}
\newcommand{\scrs}{{\mathcal{S}}}
\newcommand{\scrt}{{\mathcal{T}}}
\newcommand{\scru}{{\mathcal{U}}}
\newcommand{\scrv}{{\mathcal{V}}}
\newcommand{\scrw}{{\mathcal{W}}}
\newcommand{\scrx}{{\mathcal{X}}}
\newcommand{\scry}{{\mathcal{Y}}}
\newcommand{\scrz}{{\mathcal{Z}}}

\newcommand{\bgamma}{{\boldsymbol{\gamma}}}
\newcommand{\bsigma}{{\boldsymbol{\sigma}}}
\newcommand{\balpha}{{\boldsymbol{\alpha}}}
\newcommand{\bbeta}{{\boldsymbol{\beta}}}
\renewcommand{\pmod}[1]{\;({\rm mod}\:#1)}
\def\psibar{{\bar{\psi}}}
\def\etabar{{\bar{\eta}}}
\def\chibar{{\bar{\chi}}}
\def\xibar{{\bar{\xi}}}
\def\lambdabar{{\bar{\lambda}}}
\def\mubar{{\bar{\mu}}}
\def\varphibar{{\bar{\varphi}}}
\def\phibar{{\bar{\phi}}}
\def\cz{\overline{z}}


\newenvironment{sarray}{
	  \textfont0=\scriptfont0
	  \scriptfont0=\scriptscriptfont0
	  \textfont1=\scriptfont1
	  \scriptfont1=\scriptscriptfont1
	  \textfont2=\scriptfont2
	  \scriptfont2=\scriptscriptfont2
	  \textfont3=\scriptfont3
	  \scriptfont3=\scriptscriptfont3
	\renewcommand{\arraystretch}{0.7}
	\begin{array}{l}}{\end{array}}

\newenvironment{scarray}{
	  \textfont0=\scriptfont0
	  \scriptfont0=\scriptscriptfont0
	  \textfont1=\scriptfont1
	  \scriptfont1=\scriptscriptfont1
	  \textfont2=\scriptfont2
	  \scriptfont2=\scriptscriptfont2
	  \textfont3=\scriptfont3
	  \scriptfont3=\scriptscriptfont3
	\renewcommand{\arraystretch}{0.7}
	\begin{array}{c}}{\end{array}}

\newcommand{\bydef}{:=}
\newcommand{\defby}{=:}
\renewcommand{\implies}{\Longrightarrow}
\newcommand{\bigpartial}{\displaystyle\partial}

%
\newcommand{\ef}[1]{\, #1}     

\newcommand{\Reof}[1]{\mathfrak{Re}(#1)}
\newcommand{\Imof}[1]{\mathfrak{Im}(#1)}
\newcommand{\eval}[1]{\left\langle {#1} \right\rangle}
\newcommand{\leval}[1]{\langle {#1} \rangle}
\newcommand{\reval}[1]{\overline{#1}}

\newcommand{\sspan}{\mathrm{span}} 
\newcommand{\kker}{\mathrm{ker}}
\newcommand{\rrank}{\mathrm{rank}}

\newcommand{\bigast}[1]{\underset{#1}{\textrm{{\huge $\ast$}}}}

\newcommand{\dx}[1] {\mathrm{d}{#1}}
\newcommand{\dede}[1]{\frac{\partial}{\partial #1}}
\newcommand{\deenne}[2]{\frac{\partial^#2}{\partial #1 ^#2}}
\newcommand{\vett}[1]{#1}

\newcommand{\tinyfrac}[2] {\genfrac{}{}{}{1}{#1}{#2} }
\newcommand{\Lfrac}[2] {\genfrac{}{}{}{0}{#1}{#2} }

\newtheorem{ansatz}{Ansatz}[section]
\newtheorem{theor}{Theorem}
\newtheorem{coroll}{Corollary}

\section{Introduction}

In the theory of harmonic analysis on Euclidean Jordan algebras
(or equivalently on symmetric cones) \cite{Faraut_94},
a central role is played by the {\em Riesz distributions}\/ $\scrr_\alpha$,
which are tempered distributions
that depend analytically on a parameter $\alpha \in \C$.
One important fact about the Riesz distributions is the
necessary and sufficient condition for positivity,
due to Gindikin \cite{Gindikin_75}:

\begin{theorem} {$\!\!\!$ \bf \protect\cite[Theorem~VII.3.1]{Faraut_94} \ }
   \label{thm.riesz1}
Let $V$ be a simple Euclidean Jordan algebra of dimension $n$ and rank $r$,
with $n = r + \frac{d}{2} r(r-1)$.
Then the Riesz distribution $\scrr_\alpha$ on~$V$
is a positive measure if and only if
$\alpha = 0,\frac{d}{2},\ldots,(r-1)\frac{d}{2}$
or $\alpha > (r-1)\frac{d}{2}$.
\end{theorem}

\noindent
The ``if'' part is fairly easy, but the ``only if'' part is reputed to be deep
\cite{Gindikin_75,Faraut_94,Ishi_00}.\footnote{
   The set of values of $\alpha$ described in Theorem~\ref{thm.riesz1}
   is the so-called {\em Wallach set}\/
   \cite{Rossi_76,Wallach_79,Lassalle_87,Faraut_88,Faraut_90,Faraut_94}.
}

The purpose of this note is to give a completely elementary proof
of the ``only if'' part of Theorem~\ref{thm.riesz1},
and indeed of the following strengthening:

\begin{theorem}
   \label{thm.riesz2}
Let $V$ be a simple Euclidean Jordan algebra of dimension $n$ and rank $r$,
with $n = r + \frac{d}{2} r(r-1)$.
Then the Riesz distribution $\scrr_\alpha$ on~$V$
is a locally finite complex measure [= complex Radon measure]
if and only if $\alpha = 0,\frac{d}{2},\ldots,(r-1)\frac{d}{2}$
or $\real\alpha > (r-1)\frac{d}{2}$.
\end{theorem}

\noindent
This latter result is also essentially known
\cite[Lemma~3.3]{Hilgert_01},
but the proof given there requires some nontrivial group theory.

The idea of the proof of Theorem~\ref{thm.riesz2} is very simple:
A distribution defined on an open subset $\Omega \subset \R^n$
by a function $f \in L^1_{\rm loc}(\Omega)$
can be extended to all of $\R^n$ as a locally finite complex measure
only if the function $f$ is locally integrable
also at the boundary of $\Omega$
(Lemma~\ref{lemma.L1loc});
furthermore, this fact survives analytic continuation in a parameter
(Proposition~\ref{prop.main}).
In the case of the Riesz distribution $\scrr_\alpha$,
a simple computation using its Laplace transform
(Lemma~\ref{lemma.Delta.L1loc})
plus a bit of extra work (Lemma~\ref{lemma.uniqueness})
allows us to determine the allowed set of $\alpha$,
thereby proving Theorem~\ref{thm.riesz2}.

Theorem~\ref{thm.riesz2} thus states a necessary and sufficient
condition for $\scrr_\alpha$ to be a distribution of order 0.
It would be interesting, more generally, to determine the order of
the Riesz distribution $\scrr_\alpha$ for each $\alpha \in \C$.

It would also be interesting to know whether this approach
is powerful enough to handle the multiparameter Riesz distributions
$\scrr_{\bm{\alpha}}$ with
$\bm{\alpha} = (\alpha_1,\ldots,\alpha_r) \in \C^r$
\cite[Theorem~VII.3.2]{Faraut_94}
and/or the Riesz distributions on homogeneous cones that
are not symmetric (i.e.\ not self-dual) and hence do not arise from
a Euclidean Jordan algebra \cite{Gindikin_75,Ishi_00}.

In an Appendix I comment on a beautiful but little-known
elementary proof of Theorem~\ref{thm.riesz1} ---
which does not extend, however, to Theorem~\ref{thm.riesz2} ---
due to Shanbhag \cite{Shanbhag_88} and
Casalis and Letac \cite{Casalis_94}.

\section{A general theorem on distributions}

We assume a basic familiarity with the theory of distributions
\cite{Schwartz_66,Hormander_90}
and recall some key notations and facts.

For each open set $\Omega \subseteq \R^n$, we define
the space $\scrd(\Omega)$ of $C^\infty$ functions
having compact support in $\Omega$,
the corresponding space $\scrd'(\Omega)$ of distributions,
and the space $\scrd^{\prime k}(\Omega)$ of distributions of order $\le k$.
In particular, the space $\scrd^{\prime 0}(\Omega)$ consists
of the distributions that are given locally
(i.e.\ on every compact subset of $\Omega$)
by a finite complex measure.

Let $f \colon\, \Omega \to \C$ be a measurable function,
and extend it to all of $\R^n$ by setting $f \equiv 0$ outside $\Omega$.
We say that $f \in L^1_{\rm loc}(\Omega)$ if, for every $x \in \Omega$,
$f$ is (absolutely) integrable on some neighborhood of $x$.
Any $f \in L^1_{\rm loc}(\Omega)$ defines a distribution
$T_f \in \scrd^{\prime 0}(\Omega)$ by
\be
   T_f(\varphi)  \;=\;  \int\limits \varphi(x) \, f(x) \, dx
   \qquad \hbox{for all } \varphi \in \scrd(\Omega)
   \;.
\ee
We are interested in knowing under what circumstances 
the distribution $T_f \in \scrd^{\prime 0}(\Omega)$
can be extended to a distribution $\widetilde{T}_f \in \scrd^{\prime 0}(\R^n)$,
i.e.\ one that is locally everywhere on~$\R^n$ a finite complex measure.

\begin{lemma}
 \label{lemma.L1loc}
Let $f \colon\, \Omega \to \C$ be in $L^1_{\rm loc}(\Omega)$,
and let $T_f \in \scrd^{\prime 0}(\Omega)$ be the corresponding distribution.
Then the following are equivalent:
\begin{itemize}
   \item[(a)]  $f \in L^1_{\rm loc}(\overline{\Omega})$,
      i.e.\ for every $x \in \overline{\Omega}$,
      $f$ is integrable on some neighborhood of $x$.\footnote{
   Since this has already been assumed for $x \in \Omega$,
   the content of hypothesis (a) is that it should hold also for
   $x \in \partial\Omega$.
}
   \item[(b)]  There exists a distribution
      $\widetilde{T}_f \in \scrd^{\prime 0}(\R^n)$ that extends $T_f$
      and is supported on $\overline{\Omega}$.
   \item[(c)]  There exists a distribution
      $\widetilde{T}_f \in \scrd^{\prime 0}(\R^n)$ that extends $T_f$.
\end{itemize}
\end{lemma}

\proof
(a) $\implies$ (b):  It suffices to define
$\widetilde{T}_f(\varphi) = \int_\Omega \varphi(x) \, f(x) \, dx$
for all $\varphi \in \scrd(\R^n)$.

(b) $\implies$ (c) is trivial.

(c) $\implies$ (a):  By hypothesis, for every $x \in \partial\Omega$
and every compact neighborhood $K \ni x$,
there exists a finite complex measure $\mu_K$ supported on $K$
such that $\widetilde{T}_f(\varphi) = \int \varphi \, d\mu_K$
for every $\varphi \in \scrd(\R^n)$ with support in $K$.
But since $\widetilde{T}_f$ extends $T_f$,
the restriction of $\mu_K$ to every compact subset of $K \cap \Omega$
must coincide with the measure $f(x) \, dx$.
Since $K \cap \Omega$ is $\sigma$-compact,
this implies that
$\int\limits_{K \cap \Omega} \!\! |f(x)| \, dx
 = |\mu_K|(K \cap \Omega) < \infty$,
so that $f$ is integrable in a neighborhood of $x$.
\qed

We now extend this idea to allow for analytic dependence on a parameter.
Let $\Omega$ be an open set in $\R^n$,
let $D$ be a connected open set in $\C^m$,
and let $F \colon\, \Omega \times D \to \C$ be a continuous function
such that $F(x,\,\cdot\,)$ is analytic on $D$ for each $x \in \Omega$.
Then, for each $\lambda \in D$, define
\be
   T_\lambda(\varphi)  \;=\;
      \int\limits \varphi(x) \, F(x,\lambda) \, dx
   \qquad \hbox{for all } \varphi \in \scrd(\Omega)
   \;.
\ee

\begin{lemma}
 \label{lemma.anal}
With $F$ as above, the map $\lambda \mapsto T_\lambda$
is analytic from $D$ into $\scrd'(\Omega)$
in the sense that $\lambda \mapsto T_\lambda(\varphi)$
is analytic for all $\varphi \in \scrd(\Omega)$.
\end{lemma}

\proof
This is an immediate consequence of the hypotheses on $F$
together with standard facts about scalar-valued analytic functions
in $\C$ (either Morera's theorem or the Cauchy integral formula)
and $\C^m$ (e.g.\ the weak form of Hartogs' theorem).
\qed

{\bf Remark.}
Weak analyticity in the sense used here is actually {\em equivalent}\/
to strong analyticity:
see e.g.\ \cite[pp.~37--39, Th\'eor\`eme~1 and Remarque~1]{Grothendieck_53}
\cite[Theorems~3.1 and 3.2]{Bochnak_71}
\cite[Theorem~1]{Grosse-Erdmann_04}.
Indeed, our hypothesis on $F$ is equivalent to the even stronger statement
that the map $\lambda \mapsto F(\,\cdot\,,\lambda)$
is analytic from $D$ into the space $C^0(\Omega)$
of continuous functions on $\Omega$,
equipped with the topology of uniform convergence on compact subsets
\cite[p.~41, example (a)]{Grothendieck_53}.
But we do not need any of these facts;
weak analyticity is enough for our purposes.
\qed

\bigskip

Putting together these two lemmas, we obtain:

\begin{proposition}
 \label{prop.main}
Let $F$ be as above,
let $D_0 \subseteq D$ be a nonempty open set,
and let $\lambda \mapsto \widetilde{T}_\lambda$
be a (weakly) analytic map of $D$ into $\scrd'(\R^n)$
such that $\widetilde{T}_\lambda$ extends $T_\lambda$
for each $\lambda \in D_0$.
Then, for each $\lambda \in D$, we have:
\begin{itemize}
   \item[(a)]  $\widetilde{T}_\lambda$ extends $T_\lambda$.
   \item[(b)]  If $\widetilde{T}_\lambda \in \scrd^{\prime 0}(\R^n)$,
       then $F(\,\cdot\,,\lambda) \in L^1_{\rm loc}(\overline{\Omega})$.
\end{itemize}
\end{proposition}

\proof
(a) This is immediate by analytic continuation:
for each $\varphi \in \scrd(\Omega)$,
both $\widetilde{T}_\lambda(\varphi)$ and $T_\lambda(\varphi)$
are (by hypothesis and Lemma~\ref{lemma.anal}, respectively)
analytic functions of $\lambda$ on $D$ that coincide on $D_0$,
therefore they must coincide on all of $D$.

(b) This is immediate from (a) together with Lemma~\ref{lemma.L1loc}.
\qed

We shall apply this setup with $F(x,\lambda) = f(x)^\lambda$
where $f \colon\, \Omega \to (0,\infty)$ is a continuous function;
in fact, we shall take $f$ to be a polynomial.

\bigskip

{\bf Remark.}  Let $P$ be a polynomial that is strictly positive on $\Omega$
{\em and vanishes on $\partial\Omega$}\/,
and define for $\real \lambda > 0$ a tempered distribution
$\scrp_\Omega^\lambda \in \scrs'(\R^n)$ by the formula
\be
   \scrp_\Omega^\lambda (\varphi)
   \;=\;
   \int\limits_\Omega P(x)^\lambda \, \varphi(x) \, dx
      \qquad \hbox{for } \varphi \in \scrs(\R^n)
   \;.
 \label{def.scrp}
\ee
Then $\scrp_\Omega^\lambda$
is a tempered-distribution-valued analytic function of $\lambda$
on the right half-plane,
and it is a deep result of Atiyah, Bernstein and S.I.~Gelfand
\cite{Bernstein_69,Atiyah_70,Bernstein_72,Bjork_79}
that $\scrp_\Omega^\lambda$ can be analytically continued
to the whole complex plane as a meromorphic function of $\lambda$
with poles on a finite number of arithmetic progressions.
It is important to note that our Proposition~\ref{prop.main}
does {\em not}\/ rely on this deep result;
rather, it says that {\em whenever}\/ such an analytic continuation exists
(however it may be constructed),
the analytically-continued distribution $\scrp_\Omega^\lambda$
can be a complex measure only if
$P^\lambda \in L^1_{\rm loc}(\overline{\Omega})$.
\qed

\section{Application to Riesz distributions}

We refer to the book of Faraut and Kor\'anyi \cite{Faraut_94}
for basic facts about symmetric cones and Jordan algebras.
Let $V$ be a simple Euclidean (real) Jordan algebra
of dimension $n$ and rank $r$,
with Peirce subspaces $V_{ij}$ of dimension $d$;
recall that $n = r + \frac{d}{2} r(r-1)$.
We denote by $(x|y) = \tr(xy)$ the inner product on~$V$,
where $\tr$ is the Jordan trace and $xy$ is the Jordan product.
Let $\Omega \subset V$ be the positive cone
(i.e.\ the interior of the set of squares in $V$,
 or equivalently the set of invertible squares in $V$);
it is self-dual, i.e.\ $\Omega^* = \Omega$.
We denote by $\Delta(x) = \det(x)$ the Jordan determinant on~$V$:
it is a homogeneous polynomial of degree $r$ on~$V$,
which is strictly positive on $\Omega$ and vanishes on $\partial\Omega$,
and which satisfies \cite[Proposition~III.4.3]{Faraut_94}
\be
   \Delta(gx)  \;=\;  \Det(g)^{r/n} \, \Delta(x)
   \qquad\hbox{for $g \in G$, $x \in V$}  \;,
 \label{eq.det.quasiinv}
\ee
where $G$ denotes the identity component of the linear automorphism group
of $\Omega$ [it is a subgroup of $GL(V)$]
and $\Det$ denotes the determinant of an endomorphism.
We then have the following fundamental Laplace-transform formula:

\begin{theorem} {$\!\!\!$ \bf \protect\cite[Corollary~VII.1.3]{Faraut_94} \ }
  \label{thm.laplace}
For $y \in \Omega$ and
$\real\alpha > (r-1) \frac{d}{2} = \frac{n}{r} - 1$, we have
\be
   \int\limits_\Omega e^{-(x|y)} \, \Delta(x)^{\alpha - \frac{n}{r}} \, dx
   \;=\;
   \Gamma_\Omega(\alpha) \, \Delta(y)^{-\alpha}
 \label{eq.laplace}
\ee
where
\be
   \Gamma_\Omega(\alpha)
   \;=\;
   (2\pi)^{(n-r)/2} \,
     \prod\limits_{j=0}^{r-1} \Gamma\Bigl( \alpha - j \frac{d}{2} \Bigr)
   \;.
 \label{def.Gamma_Omega}
\ee
\end{theorem}
   
Thus, for $\real\alpha > (r-1) \frac{d}{2}$,
the function $\Delta(x)^{\alpha - \frac{n}{r}} / \Gamma_\Omega(\alpha)$
is locally integrable on $\overline{\Omega}$ and polynomially bounded,
and so defines a tempered distribution $\scrr_\alpha$ on~$V$
by the usual formula
\be
   \scrr_\alpha(\varphi)
   \;=\;
   {1 \over \Gamma_\Omega(\alpha)} \,
   \int\limits_\Omega \varphi(x) \, \Delta(x)^{\alpha - \frac{n}{r}} \, dx
      \qquad \hbox{for } \varphi \in \scrs(V)
   \;.
 \label{def.scrr}
\ee
Using \reff{eq.laplace}, a beautiful argument
--- which is a special case of Bernstein's general method
for analytically continuing distributions of the form $\scrp_\Omega^\lambda$
\cite{Bernstein_72,Bjork_79} ---
shows that the Riesz distributions $\scrr_\alpha$
can be analytically continued to the whole complex $\alpha$-plane:

\begin{theorem} {$\!\!\!$ \bf \protect\cite[Theorem~VII.2.2 et seq.]{Faraut_94}
                 \ }
The distributions $\scrr_\alpha$ can be analytically continued
to the whole complex $\alpha$-plane
as a tempered-distribution-valued entire function of $\alpha$.
Furthermore, the distributions $\scrr_\alpha$ have the following properties:
\begin{subeqnarray}
   \scrr_0  & = & \delta   \\[1mm]
   \scrr_\alpha * \scrr_\beta  & = & \scrr_{\alpha+\beta}
      \slabel{eq.riesz.convolution}   \\[1mm]
   \Delta(\partial/\partial x) \, \scrr_\alpha & = & \scrr_{\alpha-1}
      \slabel{eq.riesz.bernstein}  \\[0mm]
   \Delta(x) \, \scrr_\alpha & = &
      \!\left( \prod\limits_{j=0}^{r-1} \Bigl( \alpha - j \frac{d}{2} \Bigr)
      \!\right) \scrr_{\alpha+1}
      \slabel{eq.riesz.product}
      \label{eq.riesz}
\end{subeqnarray}
(here $\delta$ denotes the Dirac measure at 0) and
\be
   \scrr_\alpha(\varphi\circ g^{-1}) \;=\;
     \Det(g)^{\alpha r/n} \, \scrr_\alpha(\varphi)
   \qquad\hbox{for $g \in G$, $\varphi \in \scrs(V)$}
 \label{eq.riesz.quasiinv}
\ee
(in particular, $\scrr_\alpha$ is homogeneous of degree $\alpha r -n$).
Finally, the Laplace transform of $\scrr_\alpha$ is
\be
   (\scrl \scrr_\alpha)(y)  \;=\; \Delta(y)^{-\alpha}
 \label{eq.laplace.analcont}
\ee
for $y$ in the complex tube $\Omega + iV$.
\end{theorem}

\noindent
The property \reff{eq.riesz.product} is not explicitly stated in
\cite{Faraut_94}, but for $\real\alpha > (r-1) \frac{d}{2}$
it is an immediate consequence of \reff{def.Gamma_Omega}/\reff{def.scrr},
and then for other values of $\alpha$
it follows by analytic continuation
(see also \cite[Proposition~3.1(iii) and Remark~3.2]{Hilgert_01}).
Likewise, the property \reff{eq.riesz.quasiinv} is not explicitly stated in
\cite{Faraut_94}, but for $\real\alpha > (r-1) \frac{d}{2}$
it is an immediate consequence of \reff{eq.det.quasiinv}/\reff{def.scrr},
and then for other values of $\alpha$
it follows by analytic continuation
(see also \cite[Proposition~3.1(i)]{Hilgert_01}).
It follows from (\ref{eq.riesz}a,b) that the distributions $\scrr_\alpha$
are all nonzero;  and it follows from this and \reff{eq.riesz.quasiinv}
that $\scrr_\alpha \neq \scrr_\beta$ whenever $\alpha \neq \beta$.

It is fairly easy to find a {\em sufficient}\/ condition for the
Riesz distributions to be a positive measure:


\begin{proposition}
  {$\!\!\!$ \bf \protect\cite[Proposition~VII.2.3]{Faraut_94}
   (see also \protect\cite[Section~3.2]{Hilgert_01}
      \protect\cite{Lassalle_87,Bonnefoy-Casalis_90})}
 \label{prop.positive}
\begin{itemize}
   \item[(a)]  For $\alpha = k \frac{d}{2}$ with $k=0,1,\ldots,r-1$,
        the Riesz distribution $\scrr_\alpha$ is a positive measure
        that is supported on the set of elements of $\overline{\Omega}$
        of rank exactly $k$ (which is a subset of $\partial\Omega$).
   \item[(b)]  For $\alpha > (r-1) \frac{d}{2}$,
        the Riesz distribution $\scrr_\alpha$ is a positive measure
        that is supported on $\Omega$
        and given there by a density (with respect to Lebesgue measure)
        that lies in $L^1_{\rm loc}(\overline{\Omega})$.
\end{itemize}
\end{proposition}  
  
The interesting and nontrivial fact (Theorem~\ref{thm.riesz1} above)
is that the converse of Proposition~\ref{prop.positive} is also true:
the foregoing values of $\alpha$ are the {\em only}\/ ones
for which $\scrr_\alpha$ is a positive measure.
Here I shall use Proposition~\ref{prop.main}
together with the Laplace-transform formula
\reff{eq.laplace}/\reff{eq.laplace.analcont}
to provide an alternate and extremely elementary proof
of the stronger converse result stated in Theorem~\ref{thm.riesz2}.

\begin{lemma}
   \label{lemma.Delta.L1loc}
$\Delta^\lambda \in L^1_{\rm loc}(\overline{\Omega})$
if and only if $\real\lambda > -1$;
or in other words,
$\Delta^{\alpha - \frac{n}{r}} \in L^1_{\rm loc}(\overline{\Omega})$
if and only if 
$\real\alpha > (r-1) \frac{d}{2} = \frac{n}{r} - 1$.
\end{lemma}

\proof
Since $|\Delta(x)|^\lambda = \Delta(x)^{\real\lambda}$,
it suffices to consider real values of $\lambda$.

For $\lambda > - 1$ [i.e.\ $\alpha > (r-1) \frac{d}{2}$],
fix any $y \in \Omega$:  the fact that the integral \reff{eq.laplace}
is convergent, together with the fact that $x \mapsto e^{+(x|y)}$
is locally bounded, implies that
$\Delta^\lambda \in L^1_{\rm loc}(\overline{\Omega})$.

Now consider $\lambda = -1$:  again fix any $y \in \Omega$,
and let $\mu = \!\! \inf\limits_{\begin{scarray}
                                 x \in \overline{C} \\
                                 \| x \| = 1
                            \end{scarray}}
               (x|y) > 0$
where $\| \,\cdot\, \|$ is any norm on~$V$.
Choose $R > 0$ such that $|\Delta(x)| \le 1$ whenever $\|x\| \le R$.
Then
\be
   \int\limits_{\begin{scarray}
                     x \in \Omega \\
                     \|x\| \le R
                \end{scarray}}
   \!\!\!\!\!
   e^{-(x|y)} \, \Delta(x)^{-1} \, dx
   \;\;=\;\;
   \lim\limits_{\lambda \downarrow -1}
   \int\limits_{\begin{scarray}
                     x \in \Omega \\
                     \|x\| \le R
                \end{scarray}}
   \!\!\!\!\!
   e^{-(x|y)} \, \Delta(x)^{\lambda} \, dx
\ee
by the monotone convergence theorem.
We now procced to obtain a lower bound on
\be
   M_\lambda  \;\bydef
   \int\limits_{\begin{scarray}
                     x \in \Omega \\
                     \|x\| \le R
                \end{scarray}}
   \!\!\!\!\!
   e^{-(x|y)} \, \Delta(x)^{\lambda} \, dx
   \;.
\ee
For any $\beta \ge 1$, we have
\begin{subeqnarray}
   \!\!
   \int\limits_{\begin{scarray}
                     x \in \Omega \\
                     \frac{\beta}{2} R \le \|x\| \le \beta R
                \end{scarray}}
   \!\!\!\!\!
   \!\!\!\!\!
   e^{-(x|y)} \, \Delta(x)^{\lambda} \, dx
   & = &
   \beta^{n+r\lambda}
   \!\!\!\!\!
   \int\limits_{\begin{scarray}
                     x \in \Omega \\
                     \frac{R}{2} \le \|x\| \le R
                \end{scarray}}
   \!\!\!\!\!
   \!\!\!
   e^{-\beta(x|y)} \, \Delta(x)^{\lambda} \, dx
          \\[2mm]
   & \le &
   \beta^{n+r\lambda}
   e^{-(\beta-1) \frac{R}{2} \mu}
   \!\!\!\!\!
   \int\limits_{\begin{scarray}
                     x \in \Omega \\
                     \frac{R}{2} \le \|x\| \le R
                \end{scarray}}
   \!\!\!\!\!
   \!\!\!
   e^{-(x|y)} \, \Delta(x)^{\lambda} \, dx
          \\[2mm]
   & \le &
   \beta^{n+r\lambda}
   e^{-(\beta-1) \frac{R}{2} \mu}
   M_\lambda
\end{subeqnarray}
where the first equality used the homogeneity of $\Delta$.
Now sum this over $\beta = 2^k$ ($k=1,2,3,\ldots$);
the sum is convergent, and we conclude that
\be
   \int\limits_{x \in \Omega}
      e^{-(x|y)} \, \Delta(x)^{\lambda} \, dx
   \;\le\;
   C M_\lambda
\ee
for a universal constant $C < \infty$ that is independent of $\lambda$
for $-1 < \lambda \le 0$.
Since \reff{eq.laplace} tells us that
\be
   \lim\limits_{\lambda \downarrow -1}
   \int\limits_{x \in \Omega}
   \!
      e^{-(x|y)} \, \Delta(x)^{\lambda} \, dx
   \;=\;
   +\infty
\ee
due to the pole of the gamma function at $\alpha = (r-1) \frac{d}{2}$,
we conclude that $\lim\limits_{\lambda \downarrow -1} M_\lambda = +\infty$
as well.  Therefore
\be
   \int\limits_{\begin{scarray}
                     x \in \Omega \\
                     \|x\| \le R
                \end{scarray}}
   \!\!\!\!\!
   e^{-(x|y)} \, \Delta(x)^{-1} \, dx
   \;=\;
   +\infty  \;,
\ee
which proves that $\Delta^{-1} \notin L^1_{\rm loc}(\overline{\Omega})$.

Since $\Delta$ is locally bounded, it also follows that
$\Delta^\lambda \notin L^1_{\rm loc}(\overline{\Omega})$
for $\lambda < -1$.
\qed

We shall also need a uniqueness result
related to Proposition~\ref{prop.positive}(a).
If $\mu$ is a locally finite complex measure on~$V$,
we say that $\mu$ is {\em $G$-relatively invariant with exponent $\kappa$}\/
in case
\be
   \mu(gA)  \;=\; \Det(g)^\kappa \, \mu(A)
   \qquad\hbox{for $g \in G$, $A$ compact $\subseteq V$} \;.
 \label{eq.rel.inv}
\ee
In particular, every such $\mu$ is {\em $G \cap SL(V)$-invariant}\/, i.e.
\be
   \mu(gA)  \;=\; \mu(A)
   \qquad\hbox{for $g \in G \cap SL(V)$, $A$ compact $\subseteq V$} \;.
 \label{eq.SL.inv}
\ee
Now define $\Omega_k = \{x \in \overline{\Omega} \colon\, {\rm rank}(x) = k \}$,
so that $\partial\Omega = \bigcup\limits_{k=0}^{r-1} \Omega_k$
and $\Omega = \Omega_r$.
We then have the following result,
which seems to be of some interest in its own right:

\begin{lemma}
   \label{lemma.uniqueness}
\quad\par
\vspace*{-3mm}
\begin{itemize}
   \item[(a)] The group $G \cap SL(V)$ acts transitively on each
      set $\Omega_k$ $(0 \le k \le r-1)$.
   \item[(b)] Let $\mu$ be a locally finite complex measure
      that is supported on $\Omega_k$ $(0 \le k \le r-1)$
      and is $G \cap SL(V)$-invariant.
      Then $\mu$ is a multiple of $\scrr_{kd/2}$.
   \item[(c)] Let $\mu$ be a locally finite complex measure
      that is supported on $\partial\Omega$
      and is $G$-relatively invariant with some exponent $\kappa$.
      Then there exists $k \in \{0,1,\ldots,r-1\}$
      such that $\mu$ is a multiple of $\scrr_{kd/2}$
      (and hence $\kappa = kdr/2n$ if $\mu \neq 0$).
\end{itemize}
\end{lemma}

\proof
(a)  Fix a Jordan frame $c_1,\ldots,c_r$,
and let $V = \!\!\bigoplus\limits_{1 \le i \le j \le r} \!\! V_{ij}$
be the corresponding orthogonal Peirce decomposition
\cite[Theorem~IV.2.1]{Faraut_94}.
Then, for $\lambda > 0$, let
$M_\lambda = P(c_1 + \ldots + c_{r-1} + \lambda c_r) \in GL(V)$,
where $P$ is the quadratic representation \cite[p.~32]{Faraut_94}.
{}From \cite[p.~32 and Theorem~IV.2.1(ii)]{Faraut_94}
we see that $M_\lambda$ acts
as multiplication by $\lambda^2$ on the space $V_{rr}$,
as multiplication by $\lambda$ on the spaces $V_{ir}$ with $1 \le i \le r-1$,
and as the identity on the other subspaces.\footnote{
   More generally, we see that $P(\sum \lambda_i c_i)$ acts
   as multiplication by $\lambda_i \lambda_j$ on $V_{ij}$.
}
We have $M_\lambda \in G$ \cite[Proposition~III.2.2]{Faraut_94}
and $\Det(M_\lambda) = \lambda^{(r-1)d+2} = \lambda^{2n/r}$.

Now write $e_k = c_1 + \ldots + c_k$.
By construction we have $M_\lambda e_k = e_k$ for $0 \le k \le r-1$.
Now, we know \cite[Proposition~IV.3.1]{Faraut_94} that $\Omega_k = G e_k$,
so that for any $x \in \Omega_k$ there exists $g \in G$ such that $x = g e_k$.
Therefore, if we set $\lambda = \Det(g)^{-r/2n}$,
we have $x = g M_\lambda e_k$ with $g M_\lambda \in G \cap SL(V)$.

(b) follows from (a) and Proposition~\ref{prop.positive}(a)
together with a standard result about the uniqueness of invariant measures:
see e.g.\ \cite[Chapitre~7, sec.~2.6, Th\'eor\`eme~3]{Bourbaki_63},
 \cite[p.~138, Theorem~1]{Nachbin_65}
 or \cite[Theorem~7.4.1 and Corollary~7.4.2]{Wijsman_90}.

(c) is now an easy consequence, as we can write (uniquely)
$\mu = \sum_{k=0}^{r-1} \mu_k$ with $\mu_k$ supported on $\Omega_k$,
and each $\mu_k$ is $G$-relatively invariant with exponent $\kappa$
[since each set $\Omega_k$ is a separate $G$-orbit].
But in at most one case can $\kappa$ take the correct value $kdr/2n$;
so all but one of the measures $\mu_k$ must be zero.
\qed

{\bf Remarks.}
1. Assertions (a) and (b) are false when $k=r$:
the determinant $\Delta(x)$ is invariant under the action of $G \cap SL(V)$
[cf. \reff{eq.det.quasiinv}], so $G \cap SL(V)$ cannot act transitively
on $\Omega_r$;
and all the measures $\scrr_\alpha$ with $\real\alpha > (r-1) \frac{d}{2}$
are $G$-relatively invariant [hence $G \cap SL(V)$-invariant]
and supported on $\Omega_r$.

2.  A slight weakening of Lemma~\ref{lemma.uniqueness}(b)
--- in which ``$G \cap SL(V)$-invariant'' is replaced by
``$G$-relatively invariant with some exponent $\kappa$'' ---
is asserted in \cite[p.~391, Remarque~3]{Lassalle_87},
but the proof given there is insufficient
(if it were valid, it would apply also to $k=r$).
However, Michel Lassalle has kindly communicated to me
a simple alternative proof of this result,
based on \cite[Th\'eor\`eme~3 and Proposition~11(b)]{Lassalle_87}.

3.  Further information on the Riesz measures $\scrr_{kd/2}$
for $0 \le k \le r-1$ can be found in \cite{Lassalle_87,Bonnefoy-Casalis_90}.
\qed

\proofof{Theorem~\ref{thm.riesz2}}
We already know from Proposition~\ref{prop.positive}(b)
that $\scrr_\alpha$ is a locally finite complex measure
for $\real\alpha > (r-1) \frac{d}{2}$.
On the other hand,
by applying Proposition~\ref{prop.main}
to $F(x,\alpha) = \Delta(x)^{\alpha-\frac{n}{r}}/\Gamma_\Omega(\alpha)$
and using Lemma~\ref{lemma.Delta.L1loc},
we deduce that $\scrr_\alpha$ is {\em not}\/ a locally finite complex measure
whenever $\real\alpha \le (r-1) \frac{d}{2}$
and $\Gamma_\Omega(\alpha) \neq \infty$.
So it remains only to study the values of $\alpha$
for which $\Gamma_\Omega(\alpha) = \infty$,
namely $\alpha \in \{0,\frac{d}{2},\ldots,(r-1)\frac{d}{2}\} - \N$.
For $\alpha \in \{0,\frac{d}{2},\ldots,(r-1)\frac{d}{2}\}$
we know from Proposition~\ref{prop.positive}(a)
that $\scrr_\alpha$ is a positive measure.
%
%
For $\alpha \in \bigl( \{0,\frac{d}{2},\ldots,(r-1)\frac{d}{2}\} - \N \bigr)
        \setminus \{0,\frac{d}{2},\ldots,(r-1)\frac{d}{2}\}$,
we know from Proposition~\ref{prop.positive}(a) and \reff{eq.riesz.bernstein}
that $\scrr_\alpha$ is a distribution supported on $\partial\Omega$;
and by \reff{eq.riesz.quasiinv} and
Lemma~\ref{lemma.uniqueness}(b) we conclude that
it cannot be a locally finite complex measure
(here we use the fact that $\scrr_\alpha \neq \scrr_\beta$
 when $\alpha \neq \beta$).
\qed

{\bf Remark.}
For $\real\alpha < 0$, an alternate proof that $\scrr_\alpha$
is not a complex measure can be based on the following fact,
which is a special case of the $N=0$ case of \cite[Theorem~7.4.3]{Hormander_90}
(compare \cite[Theorem~7.3.1]{Hormander_90})
but can also easily be proven by direct computation:

\begin{lemma}
   \label{lemma.laplace}
Let $\Omega$ be a proper open convex cone in a real vector space $V$,
and let $\Omega^* \subset V^*$ be the open dual cone.
Let $T \in \scrs'(V) \cap \scrd^{\prime 0}(V)$
be a tempered distribution of order 0
(i.e.\ a polynomially bounded complex measure)
that is supported in $\overline{\Omega}$.
Then the Laplace transform $\scrl T$
is analytic in the complex tube $\Omega^* + iV^*$
and is bounded in every set $K + \Omega^* + iV^*$
where $K$ is a compact subset of $\Omega^*$.
\end{lemma}

\noindent
It then follows from \reff{eq.laplace.analcont}
that $\scrr_\alpha$ cannot be a locally finite complex measure
when $\real\alpha < 0$,
because $\Delta(y)^{-\alpha}$ is unbounded at infinity.
This argument handles (without the need for Lemma~\ref{lemma.uniqueness})
the cases $d=1$ (real symmetric matrices)
and $d=2$ (complex hermitian matrices)
in Theorem~\ref{thm.riesz2}.
\qed

\appendix
\section{Remarks on an elementary proof of Theorem~\ref{thm.riesz1}}

Casalis and Letac \cite[Proposition~5.1]{Casalis_94}
have given an elementary proof of Theorem~\ref{thm.riesz1}
that deserves to be more widely known than it apparently is.\footnote{
   Science Citation Index shows only ten publications citing \cite{Casalis_94},
   and six of these have an author in common with \cite{Casalis_94}.
}
They employ a method due to Shanbhag \cite[p.~279, Remark~3]{Shanbhag_88} ---
who proved Theorem~\ref{thm.riesz1} for the cases of
real symmetric and complex hermitian matrices ---
which they abstract as a general ``Shanbhag principle''
\cite[Proposition~3.1]{Casalis_94}.
Here I would like to abstract their method even further,
with the aim of revealing its utter simplicity and beauty.

Let $V$ be a finite-dimensional real vector space,
and let $V^*$ be its dual space.
We then make the following trivial observations:

(a)  If $\mu$ is a positive (i.e.\ nonnegative) measure on $V$,
then its Laplace transform
\be
   \scrl(\mu)(y)  \;=\;  \int e^{-\langle y,x \rangle} \, d\mu(x)
\ee
is nonnegative on any subset of $V^*$ where it is well-defined
(i.e.\ where the integral is convergent).

(b)  If $\mu$ is a positive measure on $V$,
then so is $f\mu$ for every continuous (or even bounded measurable)
function $f$ on $V$ that is nonnegative on $\supp \mu$.

(c) If $\mu$ is a (positive or signed) measure on $V$ whose Laplace transform
is well-defined (and finite) on a nonempty open set $\Theta \subseteq V^*$,
then the same is true for $P\mu$, where $P$ is any polynomial on $V$;
furthermore, $\scrl(P\mu) = P(-\partial) \scrl(\mu)$.\footnote{
   Indeed, the same holds if the measure $\mu$ is replaced by a
   distribution $T \in \scrd'(V)$.
   See \cite[Chapitre~VIII]{Schwartz_66}
   or \cite[Section~7.4]{Hormander_90}
   for the theory of the Laplace transform on $\scrd'(V)$.
}

Putting together these observations, we conclude:

\begin{proposition}[Shanbhag--Casalis--Letac principle]
   \label{prop.A1}
If $\mu$ is a positive measure on $V$ whose Laplace transform
is well-defined (and finite) on a nonempty open set $\Theta \subseteq V^*$,
and $P$ is a polynomial on $V$ that is nonnegative on $\supp \mu$,
then $P(-\partial) \scrl(\mu) \ge 0$ everywhere on $\Theta$.
\end{proposition}

{\bf Remark.}  Proposition~\ref{prop.A1} also has a strong converse,
which we shall state and prove at the end of this appendix.
\qed

\medskip

Using Proposition~\ref{prop.A1},
we can give the following slightly simplified version of
the Shanbhag--Casalis--Letac argument:

\bigskip\noindent
{\sc Proof of Theorem~\ref{thm.riesz1},
     based on \cite[Proposition~5.1]{Casalis_94}.\ }
In view of Proposition~\ref{prop.positive},
it suffices to prove the converse statement.
So let $\alpha \in \R$ and suppose that $\scrr_\alpha$ is a positive measure.
Using Proposition~\ref{prop.A1} with $P = \Delta$
together with the Laplace-transform formula \reff{eq.laplace.analcont},
we conclude that
\be
   \Delta(-\partial/\partial y) \, \Delta(y)^{-\alpha}  \;\ge\; 0
   \qquad\hbox{for all $y \in \Omega$} \;.
 \label{eq.proof.prop.A1.star1}
\ee
But the ``Cayley'' identity \cite[Proposition~VII.1.4]{Faraut_94}
tells us that
\be
   \Delta(\partial/\partial y) \, \Delta(y)^{\lambda}
   \;=\;
   \Delta(y)^{\lambda-1} \prod_{j=0}^{r-1} \Bigl( \lambda + j {d \over 2} \Bigr)
   \;,
 \label{eq.cayley}
\ee
hence (since $\Delta$ is homogeneous of degree $r$)
\be
   \Delta(-\partial/\partial y) \, \Delta(y)^{-\alpha}
   \;=\;
   \Delta(y)^{-\alpha-1} \prod_{j=0}^{r-1} \Bigl( \alpha - j {d \over 2} \Bigr)
   \;.
 \label{eq.proof.prop.A1.star2}
\ee
It follows from \reff{eq.proof.prop.A1.star1} and \reff{eq.proof.prop.A1.star2}
that $\scrr_\alpha$ is {\em not}\/ a positive measure when
$(r-2) {d \over 2} < \alpha <  (r-1) {d \over 2}$.
But using the convolution equation \reff{eq.riesz.convolution} with $\beta=d/2$
together with the fact that $\scrr_{d/2}$ is a positive measure
[Proposition~\ref{prop.positive}(a)],
we conclude successively that $\scrr_\alpha$ is not a positive measure
when
$(k-1) {d \over 2} < \alpha <  k {d \over 2}$
for {\em any}\/ integer $k \le r-1$.
This leaves only negative multiples of $d/2$;
and the argument given after Lemma~\ref{lemma.laplace}
shows that $\scrr_\alpha$ is not a positive measure
whenever $\alpha < 0$.\footnote{
   {\sc Alternate argument:}
   For $k=1,2,3,\ldots$ we know from Proposition~\ref{prop.positive}(a,b)
   and \reff{eq.riesz.quasiinv}
   that $\scrr_{kd/2}$ is a positive measure
   that is not supported on a single point.
   If $\scrr_{-kd/2}$ were a positive measure
   (recall that we know it is nonzero),
   then $\scrr_{kd/2} * \scrr_{-kd/2}$
   could not be supported on a single point,
   contrary to the fact that $\scrr_{kd/2} * \scrr_{-kd/2} = \delta$
   [cf.\ (\ref{eq.riesz}a,b)].
 \label{footnote_alternate}
}
\qed

{\bf Remarks.}
1. This method has been used recently by Letac and Massam
\cite[proof of Proposition~2.3]{Letac_08}
to determine the set of acceptable powers $p$
for the noncentral Wishart distribution,
generalizing the earlier proof of
Shanbhag \cite{Shanbhag_88} and Casalis and Letac \cite{Casalis_94}
for the ordinary Wishart distribution
(which is essentially Theorem~\ref{thm.riesz1}).

2.  A very different proof of Theorem~\ref{thm.riesz1}
for the cases $d=1,2$, using zonal polynomials,
was given by Peddada and Richards \cite[Theorems~1 and 3]{Peddada_91}.
\qed

\medskip

But this is not yet the end of the story;
the proof can be further simplified.
The use of the Laplace transform in the foregoing proof
is in reality a red herring,
as it is used {\em twice}\/ in opposite directions:
once in the proof of Proposition~\ref{prop.A1},
and once again in the proof of \reff{eq.cayley}.\footnote{
   The simplest proof of \reff{eq.cayley} is probably the
   one given in \cite[Proposition~VII.1.4]{Faraut_94},
   using Laplace transforms.
   However, direct combinatorial proofs are also possible:
   see \cite{CSS_cayley} for a detailed discussion
   in the cases of real symmetric and complex hermitian matrices.
}
We can therefore give a direct proof that makes almost no reference
to the Laplace transform:

\secondproofof{Theorem~\ref{thm.riesz1}}
Consider first $(r-2) {d \over 2} < \alpha <  (r-1) {d \over 2}$.
If $\scrr_\alpha$ is a positive measure,
then so is $\Delta(x) \, \scrr_\alpha$,
which by \reff{eq.riesz.product} equals $C_\alpha \scrr_{\alpha+1}$,
where
\be
   C_\alpha \;=\; \prod\limits_{j=0}^{r-1} \Bigl( \alpha - j \frac{d}{2} \Bigr)
   \;<\; 0  \;.
\ee
It follows that $\scrr_{\alpha+1}$ must be a negative (i.e.\ nonpositive)
measure.  But this is surely not the case,
as the Laplace-transform formula \reff{eq.laplace.analcont}
immediately implies that {\em no}\/ $\scrr_\beta$
can be a negative measure.\footnote{
   It would be interesting to know whether this residual use of 
   the Laplace transform can be avoided.
   For $d \le 2$ it can definitely be avoided,
   as $\alpha+1 > (r-1) {d \over 2}$,
   so that $\scrr_{\alpha+1}$ is a nonzero positive measure
   by Proposition~\ref{prop.positive}(b);
   but for $d > 2$ I do not know.
}
This shows that $\scrr_\alpha$ is not a positive measure when
$(r-2) {d \over 2} < \alpha <  (r-1) {d \over 2}$.
The proof is then completed as before.\footnote{
   The argument given after Lemma~\ref{lemma.laplace}
   explicitly uses the Laplace transform.
   But the alternate argument given in footnote~\ref{footnote_alternate}
   does not.
}
\qed

\medskip

It would be interesting to know whether this approach
is powerful enough to handle the multiparameter Riesz distributions
\cite[Theorem~VII.3.2]{Faraut_94}
and/or the Riesz distributions on homogeneous cones that
are not symmetric
and hence do not arise from
a Euclidean Jordan algebra \cite{Gindikin_75,Ishi_00}.

\medskip

To conclude, let us give the promised strong converse to
Proposition~\ref{prop.A1}:

\begin{proposition}
   \label{prop.A2}
Let $T \in \scrd'(V)$ be a distribution whose Laplace transform
is well-defined on a nonempty open set $\Theta \subseteq V^*$.
Let $S \subseteq V$ be a closed set,
and suppose that there exists $y_0 \in \Theta$ such that
$[P(-\partial) \scrl(T)](y_0) \ge 0$
for all polynomials $P$ on $V$ that are nonnegative on $S$.
Then $T$ is in fact a positive measure that is supported on $S$.
\end{proposition}

\proof
By replacing $T(x)$ by $e^{-\langle y_0,x \rangle} T(x)$,
we can assume without loss of generality that $y_0 = 0$.
Then the derivatives of $\scrl(T)$ at the origin give us the moments of~$T$;
and the hypothesis $[P(-\partial) \scrl(T)](y_0) \ge 0$ implies,
by Haviland's theorem \cite{Haviland_35,Haviland_36}
  \cite[Theorem~3.1.2]{Marshall_08},
that there exists a positive measure $\mu$ supported on $S$
that has these moments.
Furthermore, the analyticity of $\scrl(T)$ in the open set $\Theta + iV^*$
implies that these moments satisfy a bound of the form
$|c_{\bf n}| \le A B^{|{\bf n}|} {\bf n}!$,
so that $\int e^{\epsilon |x|} \, d\mu(x) < \infty$ for some $\epsilon > 0$.
It follows that the Laplace transform $\scrl(\mu)$ is well-defined
and analytic in a neighborhood of the origin;
and since its derivatives at the origin agree with those of $\scrl(T)$,
we must have $\scrl(\mu) = \scrl(T)$.
But by the injectivity of the distributional Laplace transform
\cite[p.~306, Proposition~6]{Schwartz_66},
it follows that $\mu = T$.
\qed

In Proposition~\ref{prop.A2} it is essential that the Laplace transform
of $T$ be well-defined on a nonempty open set $\Theta \ni y_0$,
or in other words (when $y_0 = 0$) that
$T$ have some exponential decay at infinity
[in the sense that
 $\cosh(\epsilon |x|) T \in \scrs'(V)$ for some $\epsilon > 0$].
It is {\em not}\/ sufficient for $T$ to have finite moments of all orders
satisfying $T(P) \ge 0$
for all polynomials $P$ on $V$ that are nonnegative on $S$.
Indeed, Stieltjes' \cite{Stieltjes_1894} famous example
\be
   f(x)  \;=\;
   \begin{cases}
     e^{-\log^2 x} \sin(2\pi \log x)  & \text{for } x > 0  \\[1mm]
     0                                & \text{for } x \le 0
   \end{cases}
\ee
belongs to $\scrs(\R)$ and has zero moments of all orders
[i.e.\ $T(P) = 0$ for all polynomials $P$]
but is not nonnegative.

\section*{Acknowledgments}

I wish to thank Jacques Faraut for an extremely helpful conversation
and correspondence,
and in particular for pointing out a gap in my original proof
of Theorem~\ref{thm.riesz2} (which used Lemma~\ref{lemma.laplace} only
and hence was valid only for the cases $d=1,2$).
I also wish to thank Michel Lassalle for correspondence
concerning Lemma~\ref{lemma.uniqueness},
Muriel Casalis for correspondence
and for providing me a copy of \cite{Bonnefoy-Casalis_90},
and Christian Berg for correspondence
concerning Proposition~\ref{prop.A2}.

I thank the Institut Henri Poincar\'e -- Centre Emile Borel
for hospitality during the programme on
Interacting Particle Systems, Statistical Mechanics and Probability Theory
(September--December 2008),
where this work was almost completed.

This research was supported in part
by U.S.\ National Science Foundation grant PHY--0424082.

\end{document}